\documentclass[11pt]{article}
\usepackage{graphicx}
\usepackage{amsmath,amsfonts,amssymb}
\usepackage{epstopdf}
\usepackage{url}
\usepackage{amsmath}
\usepackage{graphicx}
\usepackage{epsfig}
\usepackage{color}
\def\Limsup{\mathop{{\rm Lim}\,{\rm sup}}}

\def\tto{\;{\lower 1pt \hbox{$\rightarrow$}}\kern -10pt
\hbox{\raise 2pt \hbox{$\rightarrow$}}\;}
\def\Hat{\widehat}

\def\ra{\rangle}
\def\la{\langle}

\def\B{\Bbb B}
\def\h{\hfill\Box}
\def\R{\Bbb R}
\def\N{I\!\!N}
\def\ox{\bar{x}}

\def\oy{\bar{y}}
\def\oz{\bar{z}}

\def\ou{\bar{u}}

\def\co{\mbox{\rm co}\,}
\def\gph{\mbox{\rm gph}\,}
\def\epi{\mbox{\rm epi}\,}

\def\dom{\mbox{\rm dom}\,}

\def\h{\hfill\square}
\def\dn{\downarrow}

\def\ph{\varphi}

\def\st{\stackrel}

\def \N{I\!\!N}

\newcounter{lk}

\begin{document}

\begin{center}
\vspace*{0.3in} \textbf{A Unified Approach to Convex and Convexified Generalized Differentiation of Nonsmooth Functions and Set-Valued Mappings}\\[1ex]
{\em Dedicated to Professor Boris Mordukhovich on the occasion of his 65th birthday}\\[1ex]

Nguyen Mau Nam\footnote{Fariborz Maseeh Department of
Mathematics and Statistics, Portland State University, PO Box 751, Portland, OR 97207, United States (mau.nam.nguyen@pdx.edu). The research of Nguyen Mau Nam was partially supported by
the Simons Foundation under grant \#208785.}, Nguyen Dinh Hoang\footnote{Tan Tao University, Long An Province,
Vietnam (hoang1311@gmail.com)},
 R. Blake Rector\footnote{Fariborz Maseeh Department of
Mathematics and Statistics, Portland State University, PO Box 751, Portland, OR 97207, United States (r.b.rector@pdx.edu).}
\end{center}

{\small \textbf{Abstract.} In the early 1960's, Moreau and Rockafellar introduced a concept of called \emph{subgradient} for
convex functions, initiating the developments of theoretical and applied convex analysis. The needs of going beyond convexity motivated the pioneer works by Clarke considering generalized differentiation theory of Lipschitz continuous functions. Although Clarke generalized differentiation theory is applicable for nonconvex functions, convexity still plays a crucial role in Clarke subdifferential calculus. In the mid 1970's,
Mordukhovich developed another generalized differentiation theory for nonconvex functions and set-valued
mappings in which the ``umbilical cord with convexity'' no longer exists. The primary goal of this paper is to present a unified approach and shed new light on convex and Clarke generalized differentiation theories using the concepts and techniques from Mordukhovich's developments.  }

\medskip
\vspace*{0,05in} \noindent {\bf Key words.} subgradient, subdifferential, coderivative, generalized differentiation.

{\small \newtheorem{Theorem}{Theorem}[section]
\newtheorem{Proposition}[Theorem]{Proposition}
\newtheorem{Remark}[Theorem]{Remark} \newtheorem{Lemma}[Theorem]{Lemma}
\newtheorem{Corollary}[Theorem]{Corollary}
\newtheorem{Definition}[Theorem]{Definition}
\newtheorem{Example}[Theorem]{Example}
\renewcommand{\theequation}{\thesection.\arabic{equation}} }

\section{Introduction}

For centuries, differential calculus has served as an indispensable tool for science and technology, but the rise of more complex models requires new tools to deal with nondifferentiable
data. Geometric properties of convex functions and sets were the topics of study for many accomplished mathematicians such as Fenchel and Minkowski in the early 20th century. However, the beginning era of \emph{nonsmooth/variational analysis} did not start until the 1960's when Moreau and Rockafellar independently introduced a concept called \emph{subgradient} for convex functions and, together with many other mathematicians, started developing the theory of generalized differentiation for convex functions and sets.  The choice to start with convex functions and sets comes from the fact that they have several interesting properties that are important for applications to optimization. The theory is now called \emph{convex analysis}, a mature field serving as the mathematical foundation for \emph{convex optimization}.

The beauty and tremendous applications of convex analysis motivated the search for a new
theory to deal with broader classes of functions and sets where convexity is not assumed. The pioneer who
initiated this work was Clarke, a student of Rockafellar. In the early 1970's, Clarke
began to develop a generalized differentiation theory for the class
of Lipschitz continuous functions. This theory revolves around a notion called \emph{generalized gradient}, a concept that has led to numerous works, developments and applications.

The theory of generalized differentiation for nonsmooth functions relies on the interaction between the analytic and geometric properties of functions and sets. In the mid 1970's
the work of Mordukhovich brought this key idea to a very high level of generality and beauty as reflected reflected in two important objects
 of his theory of generalized differentiation: \emph{set-valued mappings} and \emph{optimal value functions}. The generalized derivative notions for nonsmooth
 functions and set-valued mappings introduced by Mordukhovich are now called under the names \emph{Mordukhovich subgradients} and \emph{coderivatives},
 respectively. Mordukhovich's generalized differentiation theory is effective for many applications, especially to
optimization theory. In spite of the nonconvexity of the Mordukhovich generalized derivative constructions, they possess well-developed
calculus rules in many important classes of Banach spaces including
reflexive spaces.

Our main goal in this paper is  to develop a
unified approach for generalized differentiatial calculus of convex
subgradients and Clarke subgradients. We show that the concepts and
techniques used by Mordukhovich are important, not only to his
generalized differentiation theory itself, but also to many other
aspects of nonsmooth analysis.

 The paper is organized as follows. In section 2, we introduce important concepts of nonsmooth analysis that are used throughout the paper. In section 3, we obtain an a subdifferential formula in the sense of convex analysis for the optimal value function under convexity, and then use it to derive many important subdifferential rules of convex analysis. Section 4 is devoted to deriving coderivative and subdifferential calculus rules in the sense of Clarke using Mordukhovich's constructions.

\section{Definitions and preliminaries}

In this section we present basic concepts and results of nonsmooth analysis to be used in next sections.

Let $X$ be a Banach space. For a function $\varphi\colon X\to
\R$ which is Lipschitz continuous around $\ox \in X$ with Lipschitz
modulus $\ell \geq 0$, the {\it Clarke generalized directional
derivative} of $\varphi$ at $\ox$ in direction $v\in X$ is defined
by
\begin{equation*}\label{gen dir der}
\varphi^{\circ}(\ox;v):=\displaystyle\limsup_{x\to \overline{x},t \dn
0}\dfrac{\varphi(x+tv)-\varphi(x)}{t}.
\end{equation*}

The generalized directional derivatives of $\varphi$ at $\ox$ are employed to define the Clarke subdifferential of the function at this point:
\begin{equation*}\label{gen gra}
\partial_{C}\varphi(\ox):=\{{x^*\in X^* \;|\; \langle x^*,v\rangle \leq
\varphi^{\circ}(\ox;v)\; \text{for any }\; v \in X}\}.
\end{equation*}

Given a nonempty closed set $\Omega$ and given a point $\ox\in \Omega$, the {\it Clarke normal cone} to $\Omega$ at $\ox$ is a subset of the dual space $X^*$ defined by
\begin{equation*}\label{cla nor}
N_C(\ox;\Omega):=\text{\rm cl}^*\{\displaystyle\bigcup_{\lambda>0}\lambda\partial_C
\text{dist}(\ox,\Omega)\},
\end{equation*}
where $\mbox{\rm dist}(\cdot; \Omega)$ is the \emph{distance function} to $\Omega$ with the representation
$$\mbox{\rm dist}(x; \Omega):=\inf\{\|x-w\|\; |\; w\in \Omega\}, x\in X.$$

With the definition of Clarke normal cones to nonempty closed sets in hand, the {\it Clarke subdifferential} $\partial_C\varphi(\ox)$ for a lower semicontinuous extended-real-valued function $\ph\colon X\to (-\infty, \infty]$ at $\ox\in \mbox{\rm dom }\ph:=\{x\in X\; |\; \ph(x)<\infty\}$  can be defined in terms of Clarke normal cones to the epigraph of the function by
\begin{equation*}
\partial_C \varphi(\ox) :=\{x^*\in X^*\; |\; (x^*,-1)\in
N_C((\overline{x},\varphi(\overline{x})); \mbox{epi}\varphi)\}.
\end{equation*}
Similarly, the {\it Clarke singular subdifferential} of $\ph$ at $\ox$ is defined by
\begin{equation*}
\partial_C^\infty \varphi(\ox) :=\{x^*\in X^*\; |\; (x^*,0)\in
N_C((\overline{x},\varphi(\overline{x})); \mbox{epi}\varphi)\}.
\end{equation*}

Another important normal cone structure of nonsmooth analysis is called the \emph{Fr\'echet normal cone} to $\Omega$ at $\ox \in \Omega$ and is defined by
\begin{equation*}\label{pre con}
\widehat N(\ox; \Omega):=\{x^* \in X^* \;|\; \displaystyle\limsup_{x
\xrightarrow{\Omega} \ox}\dfrac{\langle
x^*,x-\ox\rangle}{\|x-\ox\|}\leq 0\}.
\end{equation*}
If $\ox\notin\Omega$, we put $\widehat N(\ox; \Omega):=\emptyset$.

The {\it Mordukhovich normal cone} to $\Omega$ at $\ox$ is defined in terms of the Fr\'echet normal cone to the set around $\ox$ using the \emph{Kuratowski upper limit}:
\begin{equation*}
N(\ox; \Omega):=\Limsup_{x\to \ox}\widehat N(x; \Omega).
\end{equation*}

The \emph{Mordukhovich subdifferential} and \emph{singular subdifferential} of an extended-real-valued lower semicontinuous function $\ph\colon X\to (-\infty, \infty]$ at $\ox\in \mbox{\rm dom }\ph$ is then respectively defined by
\begin{align*}
&\partial \ph(\ox):=\{x^*\in X^*\; |\; (x^*, -1)\in N((\ox, \ph(\ox)); \epi \ph)\},\\
&\partial^\infty \ph(\ox):=\{x^*\in X^*\; |\; (x^*, 0)\in N((\ox, \ph(\ox)); \epi \ph)\}.
\end{align*}

It is well-known that when $\ph$ is a convex function, both subdifferential constructions reduce to the subdifferential in the sense of convex analysis
\begin{equation*}
\partial_C\ph(\ox)=\partial \ph(\ox)=\{x^*\in X^*\; |\; \la x^*, x-\ox\ra \leq \ph(x)-\ph(\ox)\; \mbox{\rm for all }x\in X\}.
\end{equation*}
Moreover, if the set $\Omega$ is convex, both normal cone structures reduce to the normal cone in the sense of convex analysis
\begin{equation*}
N_C(\ox; \Omega)=N(\ox; \Omega)=\{x^*\in X^*\; |\; \la x^*, x-\ox\ra \leq 0\; \mbox{\rm for all }x\in X\}.
\end{equation*}

The relation between the Mordukhovich normal cone and subdifferential constructions can also be represented by
\begin{equation*}
N(\ox;\Omega)=\partial \delta(\ox; \Omega),\; \ox\in \Omega,
\end{equation*}
where $\delta(\cdot;\Omega)$ is the indicator function associated with $\Omega$ given by $\delta(x;\Omega)=0$ if $x\in \Omega$, and $\delta(x;\Omega)=0$ otherwise. Similar relations also hold true for Clarke and Fr\'echet normal cones and subdifferentials.

Let $F\colon X\tto Y$ be a \emph{set-valued mapping} between two real
reflexive Banach spaces $X$ and $Y$. That means $F(x)$ is a subset of $Y$ for every $x\in X$. The \emph{domain} and \emph{graph} of $F$
are given respectively by
\begin{equation*}
\mbox{\rm dom }F:=\{x\in X\; |\; F(x)\neq \emptyset\} \; \mbox{\rm and }\gph F:=\{(x, y)\in X\times Y\; |\; y\in F(x)\}.
\end{equation*}
We say that $F$ has convex graph if its graph is a convex set in $X\times Y$. Similarly, we say that $F$ has closed graph if its graph is a closed set in $X\times Y$.  It is easy to see that if $B\colon X\to Y$ is a single-valued affine mapping, then $B$ has closed convex graph.

Let us continue with generalized derivative concepts for set-valued mappings introduced by Mordukhovich. The set-valued mapping $D^*F(\ox,\bar
y)\colon Y^*\rightrightarrows X^*$ defined by
\begin{equation*}
D^*F(\ox,\oy)(y^*)=\{x^*\in X^*\; |\; (x^*,-y^*)\in
N((\ox,\oy); \mbox{gph}F)\}
\end{equation*}
 is called the {\it Mordukhovich coderivative} of $F$ at $(\ox,\oy)$. By convention,
$D^*F(\ox,\oy)(y^*)=\emptyset$ for all $(\overline{x},\oy)\notin
\gph F$ and $y^*\in Y^*$. When $F$ is single-valued, one writes
$D^*F (\ox)$ instead of $D^*F(\ox,\oy)$, where $\oy=F(\ox)$. The
corresponding  {\it Clarke coderivative} is similarly defined by
\begin{equation*}
D^*_CF(\ox,\oy)(y^*)=\{x^*\in X^*\;|\; (x^*,-y^*)\in
N_C((\ox,\oy); \mbox{gph}F)\}.
\end{equation*}

In general, one has
\begin{equation*}
D^*F(\ox,\oy)(y^*)\subseteq D^*_CF(\ox,\oy)(y^*)\; \mbox{\rm for all }(\ox,\oy)\in \gph F,\; y^*\in Y^*,
\end{equation*}
where the inclusion holds as equality if $F$ has convex graph.

The following theorem which establishes the relation between Clarke and Mordukhovich generalized differentiation constructions will play an important role in our paper; see {\rm\cite{m-book1}} for the definition and proof.
\begin{Theorem}
Let $X$ be an Asplund space. The following assertions hold:

{\rm (i)} Let $\Omega \subseteq X$ be a closed set. Then

\begin{equation*}
N_C(\ox;\Omega)={\rm cl^*co}N(\ox,\Omega).
\end{equation*}

{\rm (ii)} Let $\varphi\colon X \to (-\infty, \infty]$ be lower semicontinuous
and let $\ox \in \mbox{\rm dom }\varphi$. Then
\begin{equation*}
\partial_C\varphi(\overline{x})={\rm cl^*co}[\partial\varphi(\overline{x})+\partial^\infty\varphi(\overline{x})].
\end{equation*}

In particular, if $\varphi$ is Lipschitz continuous around $\ox$,
then
\begin{equation*}
\partial_C \varphi(\ox)={\rm cl^*co} \partial \varphi(\ox).
\end{equation*}
\end{Theorem}

\section{Convex Subdifferential Calculus in Asplund Spaces via Mordukhovich Subgradients}

Our main goal in this section is to develop a unified approach for convex subdifferential calculus in Asplund spaces. We show that the concepts and techniques of nonsmooth analysis involving Mordukhovich subdifferential and coderivative constructions in Asplund spaces can be employed to shed new light on a fundamental subject of nonsmooth analysis: convex subdifferential calculus.

Let us start with some well-known examples of computing coderivatives of convex
set-valued mappings. We provide the details for the convenience of the readers. Throughout this section we assume \emph{$X,Y$ and $Z$ are Banach spaces} unless otherwise stated.

\begin{Example}\label{mr} {\rm Let $K$ be a nonempty convex set in $Y$. Define $F\colon X \tto Y$ by $F(x)=K$ for every $x\in X.$ Then
\begin{equation*}
\gph F=X \times K.
\end{equation*}
For any $(\ox,\oy)\in X \times K$,  one has $N((\ox,\oy); \gph
F)=\{0\}\times N(\oy; K)$. Thus, for $y^*\in Y^*$,
\begin{equation*}
D^*F(\ox,\oy)(y^*)=\begin{cases}
\{0\} &\text{if }\;-y^*\in N(\oy; K), \\
\emptyset & \text{otherwise}.
\end{cases}
\end{equation*}
In particular, if $K=X$, then
\begin{equation*}
D^*F(\ox,\oy)(y^*)=\begin{cases}
\{0\} &\text{if }\;y^*=0, \\
\emptyset & \text{otherwise}.
\end{cases}
\end{equation*}}
\end{Example}

The example below shows that the coderivative concept is a
generalization of the \emph{adjoint mapping} well known in
functional analysis.
\begin{Example}{\rm Let $F\colon X \tto Y$ be given by $F(x)=\{A(x)+b\}$, where $A\colon X \to Y$ is a bounded linear mapping and $b\in Y$. For $(\ox,\oy)\in \gph F$ with $\oy=A(\ox)+b$, one has
\begin{equation*}
D^*F(\ox, \oy)(y^*)=\{A^*y^*\}\; \mbox{\rm for all }y^*\in Y^*.
\end{equation*}

Let us first prove that
\begin{equation*}
N((\ox, \oy); \gph F)=\{(x^*,y^*)\in X\times Y\; |\; x^*=-A^*y^*\}.
\end{equation*}
Obviously, $\gph F$ is a convex set. By the definition, $(x^*,
y^*)\in N((\ox,\oy); \gph F)$ if and only if
\begin{equation}\label{mr1}
\la x^*, x-\ox\ra +\la y^*, F(x)-F(\ox)\ra\leq 0 \; \mbox{for all }x\in
X.
\end{equation}
One has that
\begin{align*}
\la x^*, x-\ox\ra +\la y^*, F(x)-F(\ox)\ra  &=\la x^*, x-\ox\ra +\la y^*, A(x)-A(\ox)\ra\\
&=\la x^*, x-\ox\ra +\la A^*y^*, x-\ox\ra\\
&=\la x^*+A^*y^*, x-\ox\ra.
\end{align*}
Thus, (\ref{mr1}) is equivalent to $\la x^*+A^*y^*, x-\ox\ra\leq 0$ for all $x\in X$,
which holds iff $x^*=-A^*y^*$. Now, $x^*\in D^*F(\ox, \oy)(y^*)$ if and only
if $x^*=A^*y^*$.}
\end{Example}

The example below establishes the relationship between
subdifferential and coderivative.

\begin{Example}\label{lmr}{\rm Let $f\colon X \to (-\infty, \infty]$ be a convex function. Define
\begin{equation*}
F(x):=[f(x), \infty).
\end{equation*}
Then $\gph F=\epi f$ is a convex set. For $\oy=f(\ox)\in \R$, one has
\begin{equation*}
D^*F(\ox,\oy)(\lambda)=\begin{cases}
\lambda\partial f(\ox) &\text{if }\;\lambda>0, \\
\partial^\infty f(\ox) & \text{if}\; \lambda=0,\\
\emptyset & \text{if}\; \lambda<0.
\end{cases}
\end{equation*}

First, we have the following
\begin{align*}
\gph F=\{(x,\lambda)\in X \times \R\; |\; \lambda\in
F(x)\}=\{(x,\lambda)\in X\times \R\; |\; \lambda\geq f(x)\}=\epi f.
\end{align*}
For $\lambda>0$, by the definition, $v\in D^*F(\ox,\oy)(\lambda)$ if
and only if $(v, -\lambda)\in N((\ox,\oy); \gph F)$, which is
equivalent to $(\dfrac{v}{\lambda}, -1)\in N((\ox,\oy); \epi f)$ or
$v\in \lambda\partial f(\ox)$. Similarly, for $\lambda=0$, $v\in
D^*F(\ox,\oy)(0)$ if and only if $(v, 0)\in N((\ox,\oy); \epi f)$ or
$v\in \partial^\infty f(\ox)$. Observe that if $(v, -\mu)\in N((\ox,
f(\ox)); \epi f)$, then $\mu\geq 0$. Thus, the last part of the
formula follows.}
\end{Example}

Let $F\colon X\tto Y$ be a set-valued mapping between Banach spaces, and let $\ph\colon X\times Y\to (-\infty, \infty]$ be an extended-real-valued function. The optimal value function built on $F$ and $\ph$ is given by
\begin{equation}\label{optimal value}
\mu(x):=\inf\{\ph(x,y)\; |\; y\in F(x)\}.
\end{equation}
We adopt the convention that $\inf \emptyset =\infty$. Thus, $\mu$
is an extended-real-valued function. Under convexity assumptions on
$F$ and $\ph$, we will show that convex subgradients of the optimal
value function $\mu$ can be represented by an equality in terms of
convex subgradients of the function $\ph$ and
coderivatives of the mapping $F$. Note that the result can not be derived directly from \cite{m-book1}.

The following obvious proposition guarantees that $\mu(x)>-\infty$, which is the \emph{standing assumption} for this section.
\begin{Proposition} We consider the optimal value function $\mu$ given by {\rm(\ref{optimal value})}. Suppose that $\ph(x, \cdot)$ is bounded below on $F(x)$. Then $\mu(x)>-\infty$. In particular, if there exist $b\in\R$ and a function $g\colon X\to (-\infty, \infty]$ such that for any $(x,y)\in X\times Y$ with $y\in F(x)$, one has
\begin{equation*}
\ph(x,y)\geq g(x)+b.
\end{equation*}
Then $\mu(x)>-\infty$ for all $x\in X$.
\end{Proposition}

\begin{Proposition}\label{marginal convexity} Suppose that $F$ has convex graph and $\ph$ is a convex function. Then the marginal function $\mu$ defined by {\rm (\ref{optimal value})} is a convex function.
\end{Proposition}
{\bf Proof.} Fix $x_1, x_2\in \mbox{\rm dom }\mu$ and $\lambda\in
(0,1)$. For any $\epsilon>0$, by the definition, there exist $y_i\in
F(x_i)$ such that
\begin{equation*}
\ph(x_i, y_i)<\mu(x_i) +\epsilon \; \mbox{\rm for }i=1,2.
\end{equation*}
It follows that
\begin{align*}
&\lambda\ph(x_1, y_1)<\lambda\mu(x_1) +\lambda\epsilon,\\
&(1-\lambda)\ph(x_2, y_2)<(1-\lambda)\mu(x_2)+(1-\lambda)\epsilon.
\end{align*}
Adding these inequalities and applying the convexity of $\ph$ yields
\begin{align*}
\ph(\lambda x_1+(1-\lambda) x_2, \lambda y_1+(1-\lambda)y_2)&\leq \lambda\ph(x_1,y_1)+(1-\lambda)\ph(x_2, y_2)\\
&<\lambda\mu(x_1)+(1-\lambda)\mu(x_2)+\epsilon.
\end{align*}
Since $\gph F$ is convex, $$(\lambda x_1+(1-\lambda)x_2, \lambda y_1
+(1-\lambda)y_2)=\lambda(x_1,y_1)+(1-\lambda)(x_2, y_2)\in \gph F.$$
Therefore, $\lambda y_1 +(1-\lambda)y_2\in F(\lambda
x_1+(1-\lambda)x_2)$, and hence
\begin{equation*}
\mu(\lambda x_1+(1-\lambda)x_2)\leq \ph(\lambda x_1+(1-\lambda) x_2,
\lambda
y_1+(1-\lambda)y_2)<\lambda\mu(x_1)+(1-\lambda)\mu(x_2)+\epsilon.
\end{equation*}
Letting $\epsilon\to 0$, we derive the convexity of $\mu$. $\h$

The optimal value function (\ref{optimal value}) can be used as a \emph{convex function generator} in the sense that many operations
that preserve convexity can be reduced to this function.

\begin{Proposition}\label{mr2} Consider the optimal value function {\rm (\ref{optimal value})}, where $F\colon X\tto Y$ has convex graph and $\ph\colon X\times Y\to (-\infty, \infty]$ is a convex function. Let
\begin{equation}\label{solution map}
S(\ox):=\{\oy\in F(\ox)\; |\; \mu(\ox)=\ph(\ox,\oy)\}.
\end{equation}
Assume that $\mu(\ox)<\infty$ and $S(\ox)\neq\emptyset$. For any $\oy\in S(\ox)$, one has
\begin{equation}\label{estimate 1}
\bigcup_{(u, v)\in \partial \ph(\ox,\oy)}\big[u +D^*F(\ox,\oy)(v)\big]\subseteq \partial \mu(\ox).
\end{equation}
\end{Proposition}
{\bf Proof.}  Fix any $w$ that belongs to the left side. Then there exists $(u,v)\in \partial \ph(\ox,\oy)$ such that
\begin{equation*}
w-u\in D^*F(\ox,\oy)(v).
\end{equation*}
Thus, $(w-u, -v)\in N((\ox,\oy); \gph F)$. Then
\begin{equation*}
\la w-u, x-\ox\ra -\la v, y-\oy\ra \leq 0\; \mbox{ for all }(x,y)\in
\gph F.
\end{equation*}
This implies
\begin{equation*}
\la w, x-\ox\ra \leq \la u, x-\ox\ra +\la v, y-\oy\ra \leq \ph(x,y)-\ph(\ox,\oy)=\ph(x,y)-\mu(\ox)
\end{equation*}
whenever $y\in F(x)$. It follows that
\begin{equation*}
\la w, x-\ox\ra \leq\inf_{y\in F(x)}\ph(x,y)-\mu(\ox)=\mu(x)-\mu(\ox).
\end{equation*}
Therefore, $w\in \partial \mu(\ox)$. $\h$

The opposite inclusion in (\ref{estimate 1}) also holds true under the sequentially normally compactness and the qualification condition presented below.

A nonempty closed subset $\Omega$ of a Banach space is said to be \emph{sequentially normally compact} (SNC) at $\ox\in \Omega$ if for $x_k\xrightarrow{\Omega}\ox$ and $x^*_k\in N(x_k; \Omega)$, the following implication holds:
\begin{equation*}
[x^*_k\xrightarrow{w^*}0]\Rightarrow [x^*_k\xrightarrow{\|\cdot\|}
0].
\end{equation*}
In this definition, $x_k\xrightarrow{\Omega}\ox$ means that $x_k\to \ox$ and $x_k\in \Omega$ for every $k$. Obviously, every subset of a finite dimensional Banach space is SNC. An extended-real-valued function $\ph\colon X\to (-\infty, \infty]$ is called sequentially normally epi-compact (SNEC) at $\ox$ if its epigraph is SNC at $(\ox, \ph(\ox))$.

Let us give a simple proof in the proposition below that every convex set with nonempty interior is SNC at any point of the set.
\begin{Proposition} Let $\Omega$ be a convex set with nonempty interior of a Banach space. Then $\Omega$ is SNC at any point $\ox\in \Omega$.
\end{Proposition}
{\bf Proof.} Let $\ou\in \mbox{\rm int }\Omega$ and let $\delta>0$
satisfy $\B(\ou; 2\delta)\subseteq \Omega$. Fix sequences $\{{x_k}\}$ and $\{{x^*_k}\}$ with $x_k\xrightarrow{\Omega}\ox$
and $x^*_k\in N(x_k; \Omega)$, and $x^*_k\xrightarrow{w^*}0$. Choose
$k_0$ such that $\|x_k-\ox\|<\delta$ for all $k\geq k_0$. It is not
hard to see that for any $e\in \B$, one has
\begin{equation*}
\tilde{x}_k:=\ou+\delta e +x_k -\ox \in \Omega \;\mbox{\rm for all
}k\geq k_0.
\end{equation*}
Thus,
\begin{equation*}
\la x^*_k, \tilde{x}_k-x_k \ra \leq 0.
\end{equation*}
It follows that
\begin{equation*}
\delta \la x^*_k, e\ra \leq \la x^*_k, \ox-\bar u\ra,
\end{equation*}
so $\delta \|x^*_k\|\leq \la x^*_k, \ox-\bar u\ra\to 0$. Therefore, $\|x^*_k\|\to 0$, and hence $\Omega$ is SNC at $\ox$. $\h$

\begin{Theorem}\label{mr3} Let $X$ and $Y$ be Asplund spaces. Consider the optimal value function {\rm (\ref{optimal value})}, where $F\colon X\tto Y$ has closed convex graph and $\ph\colon X\times Y\to (-\infty, \infty]$ is a lower semicontinuous convex function. Assume that $\mu(\ox)<\infty$ and $S(\ox)\neq\emptyset$, where $S$ is the solution mapping {\rm (\ref{solution map})}. For any $\oy\in S(\ox)$, one has
\begin{equation*}
\partial \mu(\ox) = \bigcup_{(u, v)\in \partial \ph(\ox,\oy)}\big[u +D^*F(\ox,\oy)(v)\big]
\end{equation*}
under the qualification condition:
\begin{equation}\label{qf}
\partial^\infty\ph(\ox,\oy)\cap [-N((\ox,\oy); \gph F)]=\{(0,0)\},
\end{equation}
and either $\ph$ is SNEC at $(\ox, \oy)$ or $\gph F$ is SNC at this point.
\end{Theorem}
{\bf Proof.} By Proposition \ref{mr2}, we only need to prove the
inclusion $\subseteq$. Fix any $w\in \partial \mu(\ox)$ and $\oy\in
S(\ox)$. Then
\begin{align*}
\la w, x-\ox\ra &\leq \mu(x)-\mu(\ox)\\
&= \mu(x)-\ph(\ox,\oy)\\
&\leq \ph(x,y)-\ph(\ox,\oy)\; \mbox{\rm for all }y\in F(x).
\end{align*}
Thus, for any $(x,y)\in X\times Y$, one has
\begin{equation*}
\la w, x-\ox\ra +\la 0, y-\oy\ra \leq \ph(x,y)+\delta((x,y); \gph F)-[\ph(\ox,\oy)+\delta((\ox,\oy); \gph F)].
\end{equation*}
Let $f(x,y):=\ph(x,y)+\delta((x,y); \gph F)$. By the definition of subdifferential, one has
\begin{equation*}
(w,0)\in  \partial f(\ox,\oy)=\partial\ph(\ox,\oy)+N((\ox,\oy); \gph F)
\end{equation*}
under the qualification condition (\ref{qf}). Thus,
\begin{equation*}
(w,0)=(u_1,v_1)+(u_2,v_2),
\end{equation*}
where $(u_1,v_1)\in \partial \ph(\ox,\oy)$ and $(u_2, v_2)\in N((\ox,\oy); \gph F)$. Then $v_2=-v_1$, and hence $(u_2, -v_1)\in N((\ox,\oy); \gph F)$. It follows that $u_2\in D^*F(\ox,\oy)(v_1)$ and
\begin{equation*}
w=u_1+u_2\in u_1+D^*F(\ox,\oy)(v_1),
\end{equation*}
where $(u_1, v_1)\in \partial \ph(\ox,\oy)$. $\h$

A function $g\colon \R^n \to (-\infty, \infty]$ is called nondecreasing componentwise if the following implication holds:
\begin{equation*}
[x_i\leq y_i \mbox{ for all } i=1,\ldots, n]\Rightarrow [g(x_1, \ldots, x_n)\leq g(y_1, \ldots, y_n)].
\end{equation*}

\begin{Proposition}\label{general com} Let $f_i\colon X\to \R$ for $i=1,\ldots,n$ be convex functions and let $h\colon X\to \R^n$ be defined by $h(x)=(f_1(x), \ldots, f_n(x))$. Suppose that $g\colon \R^n\to (-\infty, \infty]$ is a convex function that is nondecreasing componentwise. Then $g\circ h\colon X\to (-\infty, \infty]$ is a convex function.
\end{Proposition}
{\bf Proof. } Define the set-valued mapping $F\colon X\tto \R^n$ by
\begin{equation*}
F(x)=[f_1(x), \infty)\times [f_2(x), \infty)\times \ldots\times [f_n(x), \infty).
\end{equation*}
Then
\begin{equation*}
\gph F=\{(x, t_1,\ldots, t_n)\in X\times \R^n\; |\; t_i\geq f_i(x)\}.
\end{equation*}
Since $f_i$ is convex for $i=1,\ldots,n$, $\gph F$ is convex. Define $\ph\colon X\times \R^n\to (-\infty, \infty]$ by $\ph(x,y)=g(y)$. Since $g$ is nondecreasing componentwise, it is obvious that
\begin{equation*}
\inf\{ \ph(x,y)\; |\; y\in F(x)\}=g(f_1(x), \ldots, f_n(x))=(g\circ h)(x).
\end{equation*}
Thus, $g\circ h$ is convex by Proposition \ref{marginal convexity}.$\h$

\begin{Proposition} Let $X$ be an Asplund space. Let $f\colon X\to \R$ be a convex function and let $\ph\colon \R \to (-\infty, \infty]$ be a nondecreasing convex function.
Let $\ox\in X$ and let $\oy:=f(\ox)$. Then
\begin{equation*}
\partial (\ph\circ f)(\ox)=\cup_{\lambda\in \partial \ph(\oy)}\lambda \partial
f(\ox),
\end{equation*}
under the assumption that $\partial^\infty\ph(\oy)=\{0\}$ or $0\notin\partial f(\ox)$.
\end{Proposition}

{\bf Proof.}  It has been proved that $\ph\circ f$ is a convex function. Define
\begin{equation*}
F(x)=[f(x), \infty).
\end{equation*}
Since $\ph$ is a nondecreasing function, one has
\begin{equation*}
(\ph\circ f)(x)=\inf_{y\in F(x)}\ph(y).
\end{equation*}
By Theorem \ref{mr3},
\begin{equation*}
\partial (\ph\circ f)(\ox)=\cup_{\lambda \in \partial \ph(\oy)}[D^*F(\ox,\oy)(\lambda)].
\end{equation*}
Since $\ph$ is nondecreasing,  $\lambda\geq 0$ for every $\lambda\in \partial \ph(\oy)$. By Proposition \ref{lmr},
\begin{equation*}
\partial (\ph\circ f)(\ox)=\cup_{\lambda \in \partial \ph(\oy)}[D^*F(\ox,\oy)(\lambda)]=\cup_{\lambda\in \partial \ph(\oy)}\lambda \partial f(\ox).
\end{equation*}
Note that the condition $\partial^\infty\ph(\oy)=\{0\}$ or $0\notin\partial f(\ox)$ guarantees qualification condition (\ref{qf}), and $\ph$ is automatically SNCE at $\oy$. $\h$

The same approach can be applied for the general composition in Proposition \ref{general com}. A simplified version of the proposition below in finite dimensions can be found in \cite{HU}. Note that our result is new in infinite dimensions and more general in finite dimensions. Moreover, the proof is much simpler than the proof in \cite{HU}.

\begin{Proposition}\label{general com} Let $X$ be an Asplund space. Let $f_i\colon X\to \R$ for $i=1,\ldots,n$ be convex functions and let $h\colon X\to \R^n$ be defined
by $h(x)=(f_1(x), \ldots, f_n(x))$. Assume all (except possibly one) of $f_i$ for $i=1,\ldots,n$ is SNEC  at $\ox$ and the following
qualification condition holds
\begin{equation}\label{subqua1}
\big [x^*_1+\cdots +x^*_n=0,\; x^*_i\in\partial^\infty f_i(\ox)\big
]\Rightarrow \big[x^*_i=0, i=1,\ldots, n\big].
\end{equation}

Suppose that $g\colon \R^n\to (-\infty, \infty]$ is a convex function
that is nondecreasing componentwise. Then $g\circ h\colon X\to (-\infty,
\infty]$ is a convex function. Moreover, for any $\ox\in\mbox{\rm
dom }(g\circ h)$, one has
\begin{align*}
\partial (g\circ h)(\ox)=\{\sum_{i=1}^m\lambda_ix^*_i\; |\; (\lambda_1, \ldots, \lambda_m)\in \partial g(h(\ox)), x^*_i\in \partial f_i(\ox) \;\mbox{\rm for }i=1,\ldots,m\}
\end{align*}
under the condition that whenever $(\lambda_1, \ldots, \lambda_m)\in \partial^\infty g(h(\ox)), x^*_i\in \partial f_i(\ox) \;\mbox{\rm for }i=1,\ldots,m$, the following implication holds
\begin{align*}
[\sum_{i=1}^m\lambda_ix^*_i=0]\Rightarrow [\lambda_i=0\; \mbox{\rm
for }i=1,\ldots,m].
\end{align*}
\end{Proposition}
{\bf Proof.} Define the set-valued mapping $F\colon X\tto \R^n$ by
\begin{equation*}
F(x)=[f_1(x), \infty)\times [f_2(x), \infty)\times \ldots\times [f_n(x), \infty).
\end{equation*}
Then
\begin{equation*}
\gph F=\{(x, t_1,\ldots, t_n)\in X\times \R^n\; |\; t_i\geq f_i(x)\; \mbox{\rm for all }i=1,\ldots,n\}.
\end{equation*}
The set $\gph F$ is convex since since $f_i$ is convex for $i=1,\ldots,n$. Define $\ph\colon X\times \R^n\to (-\infty, \infty]$ by $\ph(x,y)=g(y)$. Since $g$ is increasing componentwise, it is obvious that
\begin{equation*}
\inf\{ \ph(x,y)\; |\; y\in F(x)\}=g(f_1(x), \ldots, f_n(x))=(g\circ h)(x).
\end{equation*}
Define
\begin{equation*}
\Omega_i=\{(x, \lambda_1, \lambda_2, \ldots, \lambda_n)\; |\; \lambda_i\geq f_i(x)\}.
\end{equation*}
Then
$$\gph F=\cap_{i=1}^n \Omega_i.$$

Using the SNC property of $\epi f_i$ for $i=1, \ldots, n$ the qualification condition
(\ref{subqua1}), and the structure of the set $\Omega_i$, one can apply
the intersection rule to get that: $(x^*, -\lambda_1, \ldots,
-\lambda_n)\in N((\ox, f_1(\ox), \ldots, f_n(\ox)); \gph F)$ if and
only if
$$(x^*, -\lambda_1, \ldots, -\lambda_n)\in N((\ox, f_1(\ox), \ldots, f_n(\ox)); \cap_{i=1}^n \Omega_i))=\sum_{i=1}^nN((\ox, f_1(\ox), \ldots, f_n(\ox)); \Omega_i).$$
If this is the case, then
\begin{equation*}
x^*=\sum_{i=1}^n x^*_i,
\end{equation*}
where $(x^*_i, -\lambda_i)\in N((\ox, f_i(\ox))$. Using Theorem \ref{mr3}, one has that $x^*\in \partial (g\circ h)(\ox)$ if and only if there exists $(\lambda_1, \ldots, \lambda_m)\in \partial g(h(\ox))$ such that
$$(x^*, -\lambda_1, \ldots, -\lambda_n)\in N((\ox, f_1(\ox), \ldots, f_n(\ox)); \gph F).$$
This is equivalent to the fact that $x^*=\sum_{i=1}^nx^*_i$, where $x^*_i\in \lambda_i\partial f_i(\ox)$. In other words,
\begin{equation*}
x^*=\sum_{i=1}^m\lambda_ix^*_i,
\end{equation*}
where $(\lambda_1, \ldots, \lambda_m)\in \partial g(h(\ox)), x^*_i\in \partial f_i(\ox) \;\mbox{\rm for }i=1,\ldots,n$. The proof is now complete. $\h$

\begin{Proposition}\label{mr5} Let $X$ and $Y$ be Asplund spaces. Let $\ph\colon X \times Y\to (-\infty, \infty]$ be an extended-real-valued function, and let $K$ be a nonempty convex subset of $Y$. Define
\begin{equation*}
\mu(x):=\inf\{\ph(x,y)\; |\; y\in K\}
\end{equation*}
and
\begin{equation*}
S(x):=\{y\in K\; |\; \ph(x,y)=\mu(x)\}.
\end{equation*}
Suppose that for every $x$, $\ph(x, \cdot)$ is bounded below on $K$.
Let $\ox\in X$ with $\mu(\ox)<\infty$ and suppose that $S(\ox)\neq
\emptyset$. Then for every $\oy\in S(\ox)$, one has
\begin{equation*}
\partial \mu(\ox)=\{u\in X^*\; |\; (u, v)\in \partial \ph(\ox,\oy), v\in N(\oy; K)\}
\end{equation*}
under the qualification condition
\begin{equation*}
[(0, v)\in \partial^\infty\ph(\ox,\oy), -v\in N(\oy; K)]\Rightarrow [v=0].
\end{equation*}
In particular, if $K=X$, then this qualification condition holds
automatically and
\begin{equation*}
\partial \mu(\ox)=\{u\in X^*\; |\; (u, 0)\in \partial \ph(\ox,\oy)\}.
\end{equation*}
\end{Proposition}
{\bf Proof.} The results follow directly from Example \ref{mr} and Theorem \ref{mr3} for the mapping $F(x)=K$.$\h$

\begin{Proposition}\label{mr5} Let $X$ and $Y$ be Asplund spaces. Let $B\colon Y\to X$ be an affine mapping with $B(y)=A(y)+b$,
where $A$ is a bounded linear mapping and $b\in Y$. Let $\ph\colon Y\to
(-\infty, \infty]$ be a convex function so that for every $x\in X$,
$\ph$ is bounded below on $B^{-1}(x)$. Define
\begin{equation*}
\mu(x):=\inf\{\ph(y)\; |\; B(y)=x\}=\inf\{\ph(y)\; |\; y\in B^{-1}(x)\}
\end{equation*}
and
\begin{equation*}
S(x):=\inf\{y\in B^{-1}(x)\;|\; \ph(y)=\mu(x)\}.
\end{equation*}
Fix $\ox\in X$ with $\mu(\ox)<\infty$ and $S(\ox)\neq\emptyset$. For
any $\oy\in S(\ox)$, one has
\begin{equation*}
\partial \mu(\ox)=(A^*)^{-1}(\partial \ph(\oy)).
\end{equation*}
\end{Proposition}
{\bf Proof.}  Let us apply Theorem \ref{mr3} for $F(x)=B^{-1}(x)$ (preimage) and $\ph(x,y)\equiv\ph(y)$. Then
\begin{equation*}
N((\ox,\oy); \gph F)=\{(u, v)\in X^*\times Y^*\; |\; -A^*u=v\}.
\end{equation*}
Thus,
\begin{equation*}
D^*F(\ox,\oy)(v)=\{u\in X^*\; |\; A^*u=v\}=(A^*)^{-1}(v).
\end{equation*}
It follows that
\begin{equation*}
\partial \mu(\ox)=\bigcup_{v\in \partial \ph(\oy)}\big [D^*F(\ox,\oy)(v)\big ]= (A^*)^{-1}(\partial \ph(\oy)).
\end{equation*}
It is not hard to verify that the qualification condition (\ref{qf}) satisfies automatically in this case.$\h$

\begin{Proposition} Let X be an Asplund space. Let $f_i\colon X \to (-\infty, \infty]$ for $i=1,2$ be convex functions. Define the convolution of $f_1$ and $f_2$ by
\begin{equation*}
(f_1\oplus f_2)(x)=\inf\{f_1(x_1)+f_2(x_2)\; |\; x_1+x_2=x\}.
\end{equation*}
Suppose that $(f_1\oplus f_2)(x)>-\infty$ for all $x\in X$ and let
$\ox\in \mbox{\rm dom }f_1\oplus f_2$. Fix $\ox_1, \ox_2\in X$ such
that $\ox=\ox_1+\ox_2$ and $(f_1\oplus f_2)(\ox)=f_1(\ox_1)+f_2(\ox_2)$.
Then
\begin{equation*}
\partial (f_1\oplus f_2)(\ox)=\partial f_1(\ox_1)\cap \partial f_2(\ox_2).
\end{equation*}
\end{Proposition}
{\bf Proof.} Let us apply Proposition \ref{mr5} for $\ph\colon X\times
X\to (-\infty, \infty]$ with $\ph(y_1, y_2)=f_1(y_1)+f_2(y_2)$ and
$A\colon X\times X\to X$ with $A(y_1, y_2)=y_1+y_2$. Then $A^*(v)=(v,v)$
for any $v\in X^*$ and $\partial \ph(\oy_1, \oy_2)=(\partial
f_1(\oy_1), \partial f_2(\oy_2))$. So $v\in
\partial (f_1\oplus f_2)(\ox)$ if and only if $A^*(v)=(v,v)\in
\partial \ph(\oy_1, \oy_2)$, i.e. $v\in \partial f_1(\ox_1)\cap
\partial f_2(\ox_2).$ $\h$

\begin{Remark}{\rm The optimal value function (\ref{optimal value}) covers many other convex operations. Thus, based on Theorem \ref{mr3}, it is possible to derive many other other convex subdifferential calculus rules. Some examples are given below.\\
{\rm (i)} Let $f_1, f_2\colon X\to (-\infty, \infty]$ be convex functions. Define
$$\ph(x,y)=f_1(x)+y$$
and $F(x)=[f_2(x), \infty)$. For any $x\in X$, one has
$$f_1(x)+f_2(x)=\inf_{y\in F(x)}\ph(x,y).$$
{\rm (ii)} Let $B\colon X\to Y$ be an affine function, and let $f\colon Y\to (-\infty, \infty]$ be a convex functions. Define $F(x)=\{B(x)\}$ and $\ph(x,y)=f(y)$ for $x\in X$ and $y\in Y$. Then
$$(f\circ B)(x)=\inf_{y\in F(x)}\ph(x,y).$$
{\rm (iii)} Let $f_1, f_2\colon X\to (-\infty, \infty]$ be convex functions. Define $F(x)=[f_1(x),\infty)\times [f_2(x), \infty)$ and
$$\ph(y_1, y_2)=\max\{y_1, y_2\}=\dfrac{|y_1-y_2|+y_1+y_2}{2}.$$
Then
$$\max\{f_1, f_2\}(x)=\inf_{(y_1, y_2)\in F(x)}\ph(y_1, y_2).$$

}
\end{Remark}

\section{Convexified Coderivative and Subdifferential Calculus}

Throughout this section, we assume that \emph{all Banach spaces under consideration are reflexive}. Under this assumption the definition of sequential normal compactness can be rewritten using weak sequential convergence. A subset $\Omega\subseteq X$ is sequentially normally compact (or shortly
SNC) at $\ox\in\Omega$ iff, for any sequences involved,  we have the implication
\begin{equation*}\label{eq:SNC}
\big[x_k\st{\Omega}{\to}\ox,\;x^*_k\st{w}{\to}0,\;x^*_k\in\Hat N(x_k;\Omega)\big]\Longrightarrow \big[\|x^*_k\|\to 0\;\mbox{ as }\;k\to\infty\big].
\end{equation*}

A subset $\Omega$ in the product space $X\times Y$ is said to be {\it
partially sequentially normally compact} (PSNC) at $(\ox,\oy)\in\Omega$
(with respect to $X$) if and only if for any sequences $(x_k,y_k) \subseteq \Omega$ and $\{{(x^*_k,y^*_k)}\}\subseteq
X^*\times Y^*$ such that $(x^*_k,y^*_k)\in \Hat N((x_k,y_k);\Omega)$,
$x^*_k\xrightarrow {w}0$ and $y^*_k\xrightarrow{\|\cdot\|} 0$, we have
$x^*_k\xrightarrow{\|\cdot\|}0$.

Accordingly, a set-valued mapping $G\colon X\tto Y$ is  SNC (PSNC) at $(\ox,\oy)\in\gph G$ iff its
graph is SNC (PSNC) at this point, and an extended-real-valued function $\ph\colon X\to(-\infty, \infty]$ is {\em
sequentially normally epi-compact} (SNEC) at $\ox\in\dom\ph$ iff its epigraph is SNC at
$(\ox,\ph(\ox))$. These properties and their partial variants are comprehensively studied and
applied in \cite{m-book1,m-book2}.

For the purposes of this paper we need the following modifications of the above properties.

\begin{Definition}\label{Def:SCNC} {\rm (i)} {\rm A set $\Omega\subseteq X$ is {\em sequentially convexly normally compact} (SCNC) at $\ox\in\Omega$ iff
we have the implication\begin{equation}\label{eq:SCNC}
\big[x_k\st{\Omega}{\to}\ox,\;x^*_k\st{w}{\to}0,\;x^*_k\in\co
N(x_k;\Omega)\big]\Longrightarrow \big[\|x^*_k\|\to 0\;\mbox{ as }\;k\to\infty\big]
\end{equation}
for any sequences involved in \eqref{eq:SCNC}. A mapping $G\colon X\tto Y$ is {\sc SCNC} at
$(\ox,\oy)\in\gph G$ iff its graph is SCNC at this point. A function $\ph\colon X\to(-\infty, \infty]$ is {\em
sequentially convexly epi-compact} (SCNEC) at $\ox\in\dom\ph$ iff its epigraph is SCNC at
$(\ox,\ph(\ox))$.

{\rm (ii)} A subset $\Omega$ of the product space $X\times Y$ is said to
be {\it partially sequentially convexifically normally compact}
(PSCNC) at $(\ox,\oy)\in\Omega$ with respect to $X$ if and only if for
any sequences $(x_k, y_k)\xrightarrow {\Omega} (\ox, \oy)$,
$\{{(x^*_k,y^*_k)}\}\subseteq X^*\times Y^*$ with
$(x^*_k,y^*_k)\in \text{co} N((x_k,y_k);\Omega)$, $x^*_k\xrightarrow
{w}0$ and $y^*_k\xrightarrow{\|\cdot\|} 0$, we have
$x^*_k\xrightarrow{\|\cdot\|}0$. A mapping $G\colon X\tto Y$ is {\sc
PSCNC} at $(\ox,\oy)\in\gph G$ iff its graph is PSCNC at this
point.}
\end{Definition}

It is easy to check that the SCNC property holds at every point of a
convex set with nonempty interior. Let us extend this result to a
broad class of nonconvex sets. Given $\Omega\subseteq X$ with $\ox\in\Omega$,
recall from \cite{r79} that $v\in X$ is a {\em hypertangent} to $\Omega$ at
$\ox$ if for some $\delta>0$ we have
\begin{equation*}\label{tan}
x+tw\in\Omega\;\mbox{ for all }\;x\in(\ox+\delta\B)\cap\Omega,\;w\in
v+\delta\B,\;\mbox{ and }\;t\in(0,\delta).
\end{equation*}

From the definition, we see that if $\Omega\subseteq X$ admits a hypertangent at $\ox$, then $\Omega$ is SCNC at this point.
Moreover, if $\ph\colon X\to (-\infty, \infty]$ is locally Lipschitz around $\ox\in \mbox{\rm dom }\ph$, then it is SCNEC at this point; see \cite{MorNamHung}
 for more details.

\begin{Lemma}\label{relation}  Let $X$ be a reflexive space and let
$\Omega$ be a subset of $X$ with $\ox\in\Omega$. Then
\begin{equation}\label{relationship}
N_C(\ox;\Omega)=\text {\rm clco} N(\ox;\Omega).
\end{equation}
\end{Lemma}
{\bf Proof.} By Theorem 3.57 in \cite{m-book1}, we have
\begin{equation*}
N_C(\ox;\Omega)=\text{cl}^*\text{co} N(\ox;\Omega).
\end{equation*}

Since in the reflexive spaces, the weak$^*$ topology and
weak topology on $X^*$ coincide, we obtain by the celebrated Mazur's
theorem that
\begin{equation*}
N_C(\ox;\Omega)=\text{cl}^*\text{co} N(\ox;\Omega)=\text{clco} N(\ox;\Omega),
\end{equation*}
which completes the proof. $\h$

Let us start with a new sum rules for Clarke coderivatives that
allow us to obtain many calculus rules for Clarke
subdifferentials and normal cones. Given $F_i\colon X\tto Y, i=1,2$, we
define a multifunction $S\colon X\times Y\tto Y\times Y$ by
\begin{equation}\label{sol}
S(x,y):=\big\{ (y_1, y_2)\in Y^{2} \; |\; y_1\in F_1(x), y_2\in
F_2(x), y_1+y_2=y\big\}.
\end{equation}

$S$ is said to be {\it inner semicontinuous} at $(\ox,\oy,
\oy_1,\oy_2)\in\gph S$ if for any sequence $\{{(x_k,y_k)}\}$ converging to $(\ox,\oy)$ with $S(x_k,y_k)\neq\emptyset$ for each $k\in \N$,
there exists $(y_{1k},y_{2k})\in S(x_k,y_k)$ such that $\{{(y_{1k},y_{2k})}\}$ contains a convergent
subsequence to $(\oy_1,\oy_2)$.

In the theorem below, we show that under the inner semicontinuity of $S$, the convexified coderivative of a sum of two set-valued mappings can be represented in terms of the coderivative of each set-valued mapping. However, this does not hold true under the so-called \emph{inner semicompactness}; see \cite{m-book1} for the definition. We use the fact that every bounded sequence in a reflexive Banach space has a subsequence that is weakly convergent, and that every closed convex set is weakly closed.

\begin{Theorem} \label{sum_rules}Let $F_i\colon X\tto Y$ for $i=1,2$ be closed-graph mappings and $(\ox,\oy)\in\gph (F_1+F_2)$. Fix $(\oy_1,\oy_2) \in S(\ox,\oy)$
such that $S$ is inner semicontinuous at $(\ox,\oy, \oy_1,\oy_2)$.
Assume that either $F_1$ is PCSNC at $(\ox,\oy_1)$ or $F_2$ is
$PCSNC$ at $(\ox,\oy_2)$, and that $\{F_1,
F_2\}$ satisfies the qualification condition
\begin{equation*}\label{qua1}
D^*_C F_1(\ox,\oy_1)(0)\cap \big (-D^*_C F_2(\ox,\oy_2)(0)\big
)=\{0\}.
\end{equation*}
Then
\begin{equation*}\label{sum}
D^*_C(F_1+F_2)(\ox,\oy)(y^*)\subseteq
D^*_CF_1(\ox,\oy_1)(y^*)+D^*_CF_2(\ox,\oy_2)(y^*).
\end{equation*}
\end{Theorem}
{\bf Proof.} Define
\begin{equation*}
\Omega_i:=\big\{ (x,y_1,y_2)\; |\; (x,y_i)\in\gph F_i\big \} \text{ for
}i=1,2.
\end{equation*}

Following the proof of \cite[Theorem 3.10]{m-book1}, one has
\begin{equation*}
\big\{ (x^*, y^*, y^*)\; | \; (x^*, y^*)\in N((\ox,\oy);\gph
(F_1+F_2))\big\} \subseteq  N((\ox,\oy_1,\oy_2);\Omega_1\cap\Omega_2).
\end{equation*}

It follows that
\begin{equation*}\label{in}
\big\{ (x^*, y^*, y^*)\; | \; (x^*, y^*)\in
\text{co}N((\ox,\oy);\gph (F_1+F_2))\big\} \subseteq
\text{co}N((\ox,\oy_1,\oy_2);\Omega_1\cap\Omega_2).
\end{equation*}

Under the assumptions made, one can apply \cite[Theorem 3.4]{m-book1} to obtain
the following inclusion:
\begin{equation*}
N((\ox,\oy_1,\oy_2);\Omega_1\cap\Omega_2) \subseteq
N((\ox,\oy_1,\oy_2);\Omega_1)+N((\ox,\oy_1,\oy_2);\Omega_2).
\end{equation*}
This implies
\begin{equation*}\label{qua2}
\text{co}N((\ox,\oy_1,\oy_2);\Omega_1\cap\Omega_2) \subseteq
\text{co}N((\ox,\oy_1,\oy_2);\Omega_1)+\text{co}N((\ox,\oy_1,\oy_2);\Omega_2).
\end{equation*}

By the definition of $\Omega_1$ and $\Omega_2$,
\begin{equation*}
N((\ox,\oy_1,\oy_2);\Omega_1)=N((\ox,\oy_1);\gph F_1)\times\{0\},
\end{equation*}
and
\begin{equation*}
N((\ox,\oy_1,\oy_2);\Omega_2)=\{(x^*, 0, y^*) \; |\; (x^*, y^*) \in
N((\ox,\oy_2);\gph F_2)\}.
\end{equation*}

Now fix any $x^* \in D^*_C(F_1+F_2)(\ox,\oy)(y^*)$. Then
\begin{equation*}
(x^*,-y^*)\in N_C((\ox,\oy); \gph(F_1+F_2)).
\end{equation*}

Let us show that
\begin{equation*}
x^* \in D^*_CF_1(\ox,\oy_1)(y^*)+D^*_CF_2(\ox,\oy_2)(y^*).
\end{equation*}

Using Lemma \ref{relation}, there exists a sequence
$\{{(x^*_k,-y^*_k)}\}$ in $\text{co}N((\ox,\oy);\gph
(F_1+F_2))$ such that $(x^*_k,-y^*_k) \to (x^*,-y^*)$, and hence $(x^*_k,-y^*_k,-y^*_k) \to (x^*,-y^*,-y^*)$.

Because $(x^*_k,-y^*_k,-y^*_k) \in
\text{co}N((\ox,\oy_1,\oy_2);\Omega_1\cap\Omega_2)$, using the above
results, there exist sequences $\{{(x^*_{1k},-y^*_k,0)}\}$ in $\text{co}N((\ox,\oy_1,\oy_2);\Omega_1)$ and
$\{{(x^*_{2k},0,-y^*_k)}\}$ in
$\text{co}N((\ox,\oy_1,\oy_2);\Omega_2)$ such that
\begin{equation*}
x^*_{1k}+x^*_{2k}=x^*_{k} \;\;\text{for every}\; k \in \N.
\end{equation*}

Without loss of generality, we assume that $F_1$ is PCSNC at
$(\ox,\oy_1)$.

We assume by a contradiction that $\{{x^*_{1k}}\}$ is
not bounded, so we can extract a subsequence, without relabeling, such
that $\|x^*_{1k}\|\rightarrow\infty$. Then
\begin{equation*}
\dfrac{1}{\|x^*_{1k}\|}(x^*_{1k}+x^*_{2k}) \rightarrow 0.
\end{equation*}
The bounded sequence $\{z^*_k\}$ defined by $z^*_k:=\dfrac{x^*_{1k}}{\|x^*_{1k}\|}$ has a weak
convergent subsequence, say, $z^*_k\xrightarrow{w} z^*$. Then
$(z^*,0,0) \in N_C((\ox,\oy_1,\oy_2);\Omega_1)$, and it is clear that
$(z^*,0,0)\in (-N_C((\ox,\oy_1,\oy_2);\Omega_2))$, which implies
$(z^*,0)\in N_C((\ox,\oy_1);\Omega_1) \bigcap (-N_C((\ox,\oy_2);\Omega_2))$.
It follows that
\begin{equation*}
z^* \in D^*_C F_1(\ox,\oy_1)(0)\cap \big (-D^*_C
F_2(\ox,\oy_2)(0)\big )=\{0\},
\end{equation*}
so $z^*=0$. Because $F_1$ is PCSNC at $(\ox,\oy_1)$, one has that
$\|z_k^*\|\rightarrow 0$. This is a contradiction since $ \left\|
{z^*_k} \right\|=1$ for every $k \in \N$.

Thus, $\{{x^*_{1k}}\}$ is bounded, so we can extract a
weak convergent subsequence. Suppose that
$x^*_{1k}\xrightarrow{w}x^*_1$ and
$x^*_{2k}\xrightarrow{w}x^*_2$.

Using the fact that
$N_C((\ox,\oy_i); \gph(F_i))=\text{cl}\text{co}N((\ox,\oy_i); \gph(F_i))$,
one obtains
\begin{equation*}
(x_i^*,-y^*) \in N_C((\ox,\oy_i); \gph(F_i)).
\end{equation*}
By definition, $x_i^* \in D^*_CF_i(\ox,\oy_i)(y^*)$  for $i=1,2$. So we
have
\begin{equation*}
x^*=x_1^*+x_2^* \in
D^*_CF_1(\ox,\oy_1)(y^*)+D^*_CF_2(\ox,\oy_2)(y^*).
\end{equation*}
The theorem has been proved. $\h$

The PCSNC holds automatically in finite dimensions, so we obtain the following sum rule for Clarke coderivatives in finite dimensions.
\begin{Corollary} \label{sum} Let $F_i\colon \R^m\tto \R^n$ for $i=1,2$ be closed-graph mappings with $(\ox,\oy)\in\gph (F_1+F_2)$. Assume
that $S$ is inner semicontinuous at $(\ox,\oy, \oy_1,\oy_2)$ and that
$\{F_1, F_2\}$ satisfies the qualification condition
\begin{equation}\label{qua1}
D^*_C F_1(\ox,\oy_1)(0)\cap \big (-D^*_C F_2(\ox,\oy_2)(0)\big
)=\{0\}.
\end{equation}
Then
\begin{equation}\label{sum}
D^*_C(F_1+F_2)(\ox,\oy)(y^*)\subseteq
D^*_CF_1(\ox,\oy_1)(y^*)+D^*_CF_2(\ox,\oy_2)(y^*).
\end{equation}
\end{Corollary}

Next we consider
\begin{equation*}
\Phi(x):= F(x)+\Delta(x;\Omega), x\in X,
\end{equation*}
where $F\colon X \tto Y$ and $\Delta(x;\Omega)=\{{0}\}\subseteq X$ if $x\in\Omega$
and $\Delta(x;\Omega)=\emptyset$ otherwise.
\begin{Proposition}\label{special} Let $\Omega$
 and $\gph F$ be closed with $\ox\in\Omega$ and $(\ox,\oy)\in\gph F$
 such that either $F$ is $PCSNC$ at $(\ox,\oy)$ or $\Omega$ is $SCNC$ at $\ox$. Assume that
 \begin{equation*}
 D^*_CF(\ox,\oy)(0)\cap (-N_C(\ox;\Omega))=\{0\}.
 \end{equation*}
 Then
 \begin{equation*}
 D^*_C(F+\Delta(\cdot;\Omega))(\ox,\oy)(y^*)\subseteq
 D^*_CF(\ox,\oy)(y^*)+N_C(\ox;\Omega), \; y^*\in Y^*.
 \end{equation*}
\end{Proposition}
{\bf Proof.} Let us apply Theorem~\ref{sum_rules} with $F_1=F$ and
$F_2=\Delta(\cdot;\Omega)$. It is not hard to see that for any $(x,y)\in
\gph (F+\Delta(\cdot;\Omega))$, we have $S(x,y)=\{(y,0)\}$, which
implies that $S$ is inner semicontinuous at $(\ox,\oy,\oy,0)$. We
also see that $\gph \Delta(\cdot;\Omega)=\Omega \times \{{0}\}$, so it is
PCSNC under the assumption that $\Omega$ is SCNC at $\ox$. We have
$N((\ox,\oy); \gph F_2)=N((\ox,\oy); \Omega \times \{{0}\})=N(\ox; \Omega)
\times Y^*$, so $N_C((\ox,\oy); \gph F_2)=N((\ox,\oy); \Omega \times
\{{0}\})=N_C(\ox; \Omega) \times Y^*$. Therefore,
$D^*_CF(\ox,\oy)(y^*)=N_C(\ox,\Omega)$. The rest follows directly from
Theorem~\ref{sum_rules}.$\h$

\begin{Corollary} Let $\Omega_1$ and $\Omega_2$ be two closed subsets $X$ and $\ox\in \Omega_1\cap\Omega_2$. Assume that $\Omega_1$ or $\Omega_2$ is $SCNC$
at $\ox$ and the qualification condition
\begin{equation*}\label{qua}
N_C(\ox;\Omega_1)\cap (-N_C(\ox;\Omega_2))=\{0\},
\end{equation*}
is satisfied. Then
\begin{equation}\label{intersect}
N_C(\ox;\Omega_1\cap\Omega_2)\subseteq N_C(\ox;\Omega_1)+N_C(\ox;\Omega_2).
\end{equation}
\end{Corollary}
{\bf Proof.} This is a special case of Proposition \ref{special} with $F_1\equiv \Delta(\cdot;\Omega_1)$.  $\h$

Let us show that the qualification condition in (\ref{qua}) is
weaker than a similar condition using Clarke tangent; see \cite{c}. Let $\Omega$ be a nonempty closed subset a Banach space $X$. For $\ox\in \Omega$, the \emph{Clarke tangent cone} $T(\ox; \Omega)$ contains all $v\in X$ such that, whenever $t_k\dn 0$ and $x_k\xrightarrow{\Omega}\ox$, there exists $w_k\to v$ with $x_k+t_kw_k\in \Omega$ for all $k$.

\begin{Proposition} Let $\Omega_1$ and $\Omega_2$ be closed subsets of $X$ and $\ox\in\Omega_1\cap\Omega_2$. Suppose that
\begin{equation}\label{tangent}
T(\ox;\Omega_1)\cap \mbox{\rm int }T(\ox;\Omega_2)\neq \emptyset.
\end{equation}
Then
\begin{equation}\label{qua}
N_C(\ox;\Omega_1)\cap (-N_C(\ox;\Omega_2))=\{0\}.
\end{equation}
\end{Proposition}
{\bf Proof.} Fix any $x^*\in N_C(\ox;\Omega_1)\cap (-N_C(\ox;\Omega_2))$.
Choose $v\in T(\ox;\Omega_1)$ such that $v+2\delta\B\subseteq
T(\ox;\Omega_2)$. Then we have $\la x^*, v\ra \leq 0$, and $\la -x^*,
v+\delta e\ra \leq 0$ for any $e\in\B$. It follows that
\begin{equation*}
\delta\la x^*, -e\ra \leq \la x^*, v\ra
\end{equation*}
for any $e\in\B$ and hence $\delta \|x^*\|\leq \la x^*, v\ra \leq
0$. Therefore $x^*=0$. $\h$

\begin{Remark}{\rm
Note that the converse of the above proposition does not hold in general. For
example, in $\R^3$, we consider $\ox=(0,0,0)$, $\Omega_1=\{{(0,0,z)\;|\;z \in
\R}\}$ and $\Omega_2=\{{(x,y,0)\;|\;x,y\in \R}\}$. We have
$N_C(\ox;\Omega_1)=\{{(x,y,0)\;|\;x,y\in \R}\}$ and
$N_C(\ox;\Omega_2)=\{{(0,0,z)\;|\;z \in \R}\}$. Thus,
\begin{equation}
N_C(\ox;\Omega_1)\cap (-N_C(\ox;\Omega_2))=\{0\},
\end{equation}
but $T(\ox;\Omega_1)\cap \text{\rm int }T(\ox;\Omega_2)= \emptyset$ since
$\text{\rm int }T(\ox;\Omega_2)=\text{\rm int }T(\ox;\Omega_1)=\emptyset$.
Therefore, we obtain a stronger finite-dimensional version of
Corollary 2.9.8 in {\cite{c}}. }
\end{Remark}

\begin{Definition}{\rm  We say that an
extended-real-valued function $\ph\colon X\to (-\infty, \infty]$ is \emph{lower regular} at $\ox\in\mbox{\rm dom }\ph$ if
$\Hat\partial\ph(\ox)=\partial_C\ph(\ox)$, where $\Hat\partial\ph(\ox)$ is the \emph{Fr\'echet normal subdifferential} of $\ph$ at $\ox$ defined by
\begin{equation*}
\Hat\partial \ph(\ox):=\big\{x^*\in X^*\;\big |\; \liminf_{x\to \ox}\dfrac{\ph(x)-\ph(\ox)-\la x^*, x-\ox\ra}{\|x-\ox\|}\geq 0\big\}.
\end{equation*}}
\end{Definition}
It is clear that any convex function is lower regular.

Now we apply the results obtained from coderivative calculus in
Theorem \ref{sum_rules} to obtain calculus for Clarke subdifferential in
reflexive Banach spaces.
\begin{Theorem}\label{subsum} Let $\ph_i\colon X\rightarrow(-\infty, \infty]$ for $i=1,\ldots,m$ be l.s.c.
around $\ox$ and finite at this point. Assume all (except possibly one) of
$\ph_i, i=1,\ldots,m$ are $SCNEC$ at $(\ox,\ph_i(\ox))$ and
\begin{equation}\label{subqua}
\big [x^*_1+\cdots+x^*_m=0,\; x^*_i\in\partial^\infty_C\ph_i(\ox)\big
]\Rightarrow x^*_i=0, i=1,\ldots, m.
\end{equation}
Then
\begin{equation}\label{sum1}
\partial_C(\ph_1+\cdots+\ph_m)(\ox)\subseteq\partial_C\ph_1(\ox)+\cdots+\partial_C\ph_k(\ox),
\end{equation}
and
\begin{equation}\label{sum2}
\partial^\infty_C(\ph_1+\cdots+\ph_m)(\ox)\subseteq\partial^\infty_C\ph_1(\ox)+\cdots+\partial^\infty_C\ph_m(\ox).
\end{equation}
The equality in {\rm(\ref{sum1})} holds if all $\ph_i$ for $i=1,\ldots,m$ are upper regular
at $\ox$.
\end{Theorem}

{\bf Proof.} We first consider the case where  $m=2$. Let us consider
$F_1(x)=[\ph_1(x),\infty)$ and $F_2(x)=[\ph_2(x),\infty)$. Obviously, $\gph F_i=\epi(\ph_i)$ for $i=1,2$ and
$\gph(F_1+F_2)=\epi(\ph_1+\ph_2)$. We see that at the point
$(\ox,\oy)$ where $\oy=\ph_1(\ox)+\ph_2(\ox)$, we have
$S(\ox,\oy)=\{(\ph_1(\ox),\ph_2(\ox))\}$, and $S$ is inner
semicontinuous at $(\ox,\oy,\ph_1(\ox),\ph_2(\ox))$. Indeed, for
every sequence $\{{(x_k,y_k)}\}$ converging to $(\ox,\oy)$ with $S(x_k,y_k)\neq\emptyset$. Fix $(\lambda_{1k}, \lambda_{2k})\in S(x_k, y_k)$. Then
\begin{equation*}
\lambda_{1k}\geq \ph_1(x_k), \lambda_{2k}\geq \ph_2(x_k), \lambda_{1k}+\lambda_{2k}=y_k.
\end{equation*}
Since $\ph_1$ and $\ph_2$ are lower semicontinous,
\begin{equation*}
\liminf_{k\to \infty}\lambda_{ik}\geq \liminf_{k\to \infty}\ph_i(x_k)\geq \ph_i(\ox)\; \mbox{\rm for }i=1,2.
\end{equation*}
We can see that $\{\lambda_{1k}\}$ and $\{\lambda_{2k}\}$ are bounded sequences. Let $\lambda_1:=\liminf_{k\to\infty} \lambda_{1k}$ and $\lambda_2:=\liminf_{k\to\infty} \lambda_{2k}$. Then there exist subsequences of $\{\lambda_{1k}\}$ and $\{\lambda_{2k}\}$ that converge to $\lambda_1$ and $\lambda_2$, respectively. Since $\lambda_1\geq \ph_1(\ox)$, $\lambda_2\geq \ph_2(\ox)$, and $\lambda_1+\lambda_2=\oy$, we see that $\lambda_1=\ph_1(\ox)$ and $\lambda_2=\ph_2(\ox)$, so $S$ is inner semicontinuous at $(\ox,\oy,\ph_1(\ox),\ph_2(\ox))$.

Because $D^*_C
F_i(\ox,\ph_i(\ox))(0)=\partial_C^\infty\ph_i(\ox)$, the qualification
condition (\ref{qua1}) holds. Applying Theorem \ref{sum}, we have
\begin{equation*}
D^*_C(F_1+F_2)(\ox,\oy)(1)\subseteq
D^*_CF_1(\ox,\ph_1(\ox))(1)+D^*_CF_2(\ox,\ph_2(\ox))(1),
\end{equation*}
and
\begin{equation*}
D^*_C(F_1+F_2)(\ox,\oy)(0)\subseteq
D^*_CF_1(\ox,\ph_1(\ox))(0)+D^*_CF_2(\ox,\ph_2(\ox))(0)
\end{equation*}
which imply (\ref{sum1}) and (\ref{sum2}) for $m=2$.

Let us prove (\ref{sum1}) holds as an equality under the regularity
conditions. Assume that all $\ph_i$ for $i=1,2$ are lower regular at
$\ox$. Then
\begin{align*}
\partial_C(\ph_1+\ph_2)&\subseteq
\partial_C\ph_1(\ox)+\partial_C\ph_2(\ox)\\
&=\Hat\partial\ph_1(\ox)+\Hat\partial\ph_2(\ox)\\
&\subseteq \Hat\partial(\ph_1+\ph_2)(\ox)\\
&\subseteq \partial_C(\ph_1+\ph_2)(\ox),
\end{align*}
which implies the equality. The proof for $m>2$ follows easily
by induction. $\h$

\begin{Proposition}
Let $\ph_i\colon X\to \R$ satisfy
\begin{equation*}
\{ v \in X\; |\;  \ph^\circ_1(\ox; v)<\infty\} \cap \mbox{\rm int} \{ v\in X\; |\;
\ph^\circ_2(\ox;v)<\infty\}\neq\emptyset.
\end{equation*}
Then
\begin{equation*}
\partial_C^\infty \ph_1(\ox)\cap (-\partial_C^\infty \ph_2(\ox))=\{0\}.
\end{equation*}
\end{Proposition}
{\bf Proof.} Fix any $x^*\in \partial_C^\infty \ph_1(\ox)\cap
(-\partial_C^\infty \ph_2(\ox))$. Then $(x^*, 0)\in N_C((\ox,
\ph_1(\ox));\epi \ph_1)$ and $(-x^*, 0)\in N_C((\ox,
\ph_2(\ox));\epi \ph_2)$. Fix an element $v\in X$ such that
$\ph_1^\circ(\ox;v)<\infty$, while $v\in \mbox{\rm int }\{ v_1\in
X\; |\; {\ph_2}^\circ(\ox;v_1)<\infty\}\neq\emptyset$. Let us show
that $x^*=0$. From the proof of \cite[Theorem 2.9.5]{c}, one finds
$\beta\in \R$ such that $(v,\beta)\in \mbox{\rm int } T((\ox,
\ph_2(\ox)); \epi \ph_2)$. From \cite[Theorem 2.9.1]{c}, we see that
if  $\gamma\in\R$ is fixed such that $\gamma > \ph_1^\circ
 (\ox;v)$, then $(v, \gamma)\in T((\ox, \ph_1(\ox));\epi \ph_1)$. Thus, there exists $\delta>0$ with
 \begin{equation*}
 (v,\beta)+\delta\B\subseteq T((\ox, \ph_2(\ox));\epi \ph_2),
 (v,\gamma)\in T((\ox, \ph_1(\ox));\epi \ph_1).
 \end{equation*}
 Then
 \begin{equation*}
 \la (x^*,0), (v,\beta)+\delta/2(e_1,e_2)\ra \leq 0\leq \la (x^*,0),
 (v,\gamma)\ra
 \end{equation*}
 whenever $\|(e_1,e_2)\|\leq 1$. It follows that $\la x^*,\delta/2
 e_1\ra \leq 0$ whenever $\|e_1\|\leq 1$, which implies $x^*=0$.
 $\h$

 The proposition below shows that the SCNEC condition for extended-real-valued functions holds under the \emph{directional Lipschitz condition}; see \cite[Definition 2.9.2]{c}.

 \begin{Proposition} Assume
 that $\ph\colon X\to (-\infty, \infty]$ is directionally Lipschitz at $\ox\in\mbox{\rm dom }\ph$. Then $\ph$ is SCNEC  at
 $(\ox, \ph(\ox))$.
 \end{Proposition}
{\bf Proof.} Assume that $\ph$ is directionally Lipschitz at $\ox$.
Then by \cite[Proposition 2.9.3]{c}, there exists $\beta\in \R$ such
that $(v,\beta)$ is a hypertangent to $\epi \ph$ at $(\ox,
\ph(\ox))$. Therefore, $\epi \ph$ is SCNC at this point. $\h$

\begin{Theorem}\label{chain1}
Let $G\colon X\tto Y,\; F\colon Y\tto Z $ be closed-graph mappings with $\oz\in
(F\circ G)(\ox)$, and
\begin{equation*}
S(x,z):=G(x)\cap F^{-1}(z)=\{{y\in G(x)|z\in F(y)}\}.
\end{equation*}
Given $\oy\in S(\ox,\oz)$, assume that $S$ is inner semicontinuous
at $(\ox,\oz,\oy)$, that either $F$ is PSCNC  at
$(\oy,\oz)$ or $G$ is PSCNC at $(\ox,\oy)$, and that
the qualification condition
\begin{equation*}
D^*_CF(\oy,\oz)(0)\cap \mbox{\rm ker }D_C^*G(\ox,\oy)=\{0\},
\end{equation*}
is fulfilled. Then
\begin{equation}\label{chain}
D^*_C(F\circ G)(\ox,\oz)(z^*)\subseteq D^*_CG(\ox,\oy)\circ
D^*_CF(\oy,\oz)(z^*),
\end{equation}
for any $z^*\in Z^*$.
\end{Theorem}
{\bf Proof.} Consider the set-valued mapping $\Phi\colon X \times Y \tto
Z$ as follows
\begin{equation*}
\Phi(x,y)=F(y)+\Delta((x,y);\gph G).
\end{equation*}
Using \cite[Theorem 1.64]{m-book1}, because $S$ is inner
semicontinuous at $(\ox,\oz,\oy)$, we have
\begin{equation*}
D^*(F\circ G)(\ox,\oz)(z^*)\subseteq \big\{ x^*\in X^*\; |\;
(x^*,0)\in D^*\Phi(\ox,\oy,\oz)(z^*)\}.
\end{equation*}
Thus,
\begin{equation*}
\mbox{\rm co }D^*(F\circ G)(\ox,\oz)(z^*)\subseteq \big\{ x^*\in X^*\; |\;
(x^*,0)\in \mbox{\rm co }D^*\Phi(\ox,\oy,\oz)(z^*)\}.
\end{equation*}
Therefore,
\begin{equation*} D^*_C(F\circ G)(\ox,\oz)(z^*)\subseteq \big\{
x^*\in X^*\; |\; (x^*,0)\in D^*_C\Phi(\ox,\oy,\oz)(z^*)\}.
\end{equation*}
It follows from Proposition~\ref{special} that
\begin{equation*}
D^*_C\Phi(\ox,\oy,\oz)(z^*) \subseteq
D^*_CF(\oy,\oz)(z^*)+N_C((\ox,\oy); \gph G).
\end{equation*}
Take $x^* \in D^*_C(F\circ G)(\ox,\oz)(z^*)$. Then $(x^*,0) \in
D^*_C\Phi(\ox,\oy,\oz)(z^*)$. Since $F=F(y)$, there exists
$y^* \in Y^*$ such that $(y^*,-z^*) \in N_C((\oy,\oz); \gph F)$ and
$(x^*,-y^*) \in N_C((\ox,\oy);\gph G)$. Then $y^* \in
D^*_CF(\oy,\oz)(z^*)$ and $x^* \in D^*_CG(\ox,\oy)(y^*)$. Therefore,
$x^*  \in D_C^* G(\ox ,\oy ) \circ D_C^*
F(\oy ,\oz )(z^* )$. The theorem has been proved. $\h$

\begin{Corollary}\label{inverse} Let $F\colon
X\tto Y$ be a closed-graph mapping, and let $\Omega\subseteq Y$ be a closed set. For $(\ox,\oy) \in \gph F$
and $\oy \in \Omega$, define
\begin{equation*}
F^{-1}(\Omega):=\big\{ x\in X\; |\; F(x)\cap\Omega\neq\emptyset\big\}.
\end{equation*}
Assume that $S(x):=F(x)\cap\Omega$ is inner semicontinuous at
$(\ox,\oy)$, and that either
$F$ is PSCNC at $(\ox,\oy)$ or $\Omega$ is SCNC
at $\oy$. Under the qualification condition
\begin{equation*}
N_C(\oy;\Omega)\cap \mbox{\rm ker }D^*_CF(\ox,\oy)=\{0\},
\end{equation*}
one has
\begin{equation*}
N_C(\ox; F^{-1}(\Omega))\subseteq D^*_CF(\ox,\oy)(N(\oy;\Omega)).
\end{equation*}
\end{Corollary}
{\bf Proof.} This follows from Theorem~\ref{chain1} with $F_1\colon Y\tto
\R$ and $G_1\colon X \tto Y$, where $G_1=F$ and
$F_1(\cdot)=\Delta(\cdot;\Omega)$. Then $\gph F_1=\Omega \times \{{0}\}$.
Obviously, if $\Omega$ is SCNC  at $\oy$, then $\gph F_1$ is PSCNC at
$(\oy,0)$. Applying Theorem~\ref{chain1}, we obtain the result.$\h$

\begin{Theorem}\label{composition} Let $f:=g\circ F$, where $F\colon X\to Y$ is a strictly
differentiable mapping and $g\colon Y\to (-\infty, \infty]$ is an
extended-real-valued function. Assume that $g$ is l.s.c
around $F(\ox)$, and that
\begin{equation*}
\partial^{\infty}_Cg(F(\ox))\cap (\mbox{\rm ker }\nabla F(\ox)^*)=\{0\}.
\end{equation*}
Then
\begin{equation*}
\partial_C f(\ox)\subseteq \nabla F(\ox)^*(\partial_C g(F(\ox))),
\end{equation*}
and
\begin{equation*}
\partial_C^\infty f(\ox)\subseteq \nabla F(\ox)^*(\partial_C^\infty g(F(\ox))).
\end{equation*}
The first inclusion holds as an equality if $g$ is lower regular at
$F(\ox)$.
\end{Theorem}
{\bf Proof.} Let us consider $E_g(x)=[g(x),\infty)$. Obviously,
\begin{equation*}
E_f(x)=(E_g\circ F)(x).
\end{equation*}
Observe that $D^*_CF(\ox)=\nabla F(\ox)^*$. The result then follows
directly from Theorem \ref{chain1} with $z^*=1$ and $z^*=0$. $\h$

\end{document}